\documentclass[12pt]{article}
\usepackage{amssymb}
\usepackage{amsfonts}
\usepackage{amsmath}
\usepackage{amsbsy}
\newcommand{\be}{\begin{equation}}
\newcommand{\ee}{\end{equation}}
\newcommand{\bea}{\begin{eqnarray}}
\newcommand{\eea}{\end{eqnarray}}
\newcommand{\ba}{\begin{array}}
\newcommand{\ea}{\end{array}}

\newcommand{\bc}{\begin{center}}
\newcommand{\ec}{\end{center}}
\newcommand{\ben}{\begin{enumerate}}
\newcommand{\een}{\end{enumerate}}
\newcommand{\bfi}{\begin{figure}}
\newcommand{\efi}{\end{figure}}

\newcommand{\bq}{\begin{quote}}
\newcommand{\eq}{\end{quote}}
\newcommand{\bqu}{\begin{quotation}}
\newcommand{\equ}{\end{quotation}}
\newenvironment{emphit}{\begin{itemize}}{\end{itemize}}
\newcommand{\bemp}{\begin{emphit}}
\newcommand{\eemp}{\end{emphit}}

\newcommand{\bt}{\begin{tabular}}
\newcommand{\et}{\end{tabular}}

\newtheorem{myth}{Theorem}[section]
\newtheorem{mylem}{Lemma}[section]
\newtheorem{mycor}{Corollary}[section]

\newtheorem{mydef}{Definition}[section]
\newtheorem{myrem}{Remark}[section]

\begin{document}
\date{}
\title{On the Derivatives of Central Loops
}
\author{T. G. Jaiy\'e\d ol\'a\\
Department of Mathematics,\\
Obafemi Awolowo University,\\
Ile Ife 220005, Nigeria.\\
jaiyeolatemitope@yahoo.com\\tjayeola@oauife.edu.ng \and
J. O. Ad\'en\'iran\thanks{All correspondence to be addressed to this author.}  \\
Department of Mathematics,\\
University of Abeokuta, \\
Abeokuta 110101, Nigeria.\\
ekenedilichineke@yahoo.com\\
adeniranoj@unaab.edu.ng} \maketitle

\begin{abstract}
The right(left) derivative, $a^{-1},e-$ and $e,a^{-1}-$ isotopes of
a C-loop are shown to be C-loops. Furthermore, for a central loop
$(L,F)$, it is shown that $\big\{F,F^{a^{-1}},F_{a^{-1},e}\big\}$
and $\big\{F,F_{a^{-1}},F_{e,a^{-1}}\big\}$ are systems of isotopic
C-loops that obey a form of generalized distributive law. Quasigroup
isotopes $(L,\otimes )$ and $(L,\ominus )$ of a loop $(L,\theta )$
and its parastrophe $(L,\theta ^*)$ respectively are proved to be
isotopic if either $(L,\otimes )$ or $(L,\ominus )$ is commutative.
If $(L,\theta )$ is a C-loop, then it is shown that $\big\{(L,\theta
),(L,\theta ^*),(L,\otimes ),(L,\oplus )\big\}$ is a system of
isotopic C-quasigroup under the above mentioned condition. It is
shown that C-loops are isotopic to some finite indecomposable groups
of the classes ${\cal D}_i,i=1,2,3,4,5$ and that the center of such
C-loops have a rank of 1,2 or 3.
\end{abstract}

\section{Introduction}
\paragraph{}Central loops('C-loops' for short) are the least studied of the Bol-Moufang type loops.
Closely related to them are left central loops(LC-loops) and right
central loops(RC-loops)(\cite{fen2}). They have been found to
behave in a manner similar to that of Moufang loops. The
identities and the equivalent ones that describe these
loops(LC,RC,C) are available in \cite{fen2} and we give them
below.
\begin{equation}\label{c:def}
(yx\cdot x)z=y(x\cdot xz)~~\textrm{central identity}
\end{equation}
\begin{equation}\label{lc1:def}
xx\cdot yz=(x\cdot xy)z~~\textrm{left central identity}
\end{equation}
\begin{equation}\label{lc2:def}
(x\cdot xy)z=x(x\cdot yz)~~\textrm{left central identity}
\end{equation}
\begin{equation}\label{lc3:def}
(xx\cdot y)z=x(x\cdot yz)~~\textrm{left central identity}
\end{equation}
\begin{equation}\label{rc1:def}
yz\cdot xx=y(zx\cdot x)~~\textrm{right central identity}
\end{equation}
\begin{equation}\label{rc2:def}
(yz\cdot x)x=y(zx\cdot x)~~\textrm{right central identity}
\end{equation}
\begin{equation}\label{rc3:def}
(yz\cdot x)x=y(z\cdot xx)~~\textrm{right central identity}
\end{equation}
They have been studied by Phillips and Vojt\v echovsk\'y
\cite{phi}, \cite{phi1}, \cite{phi2}, Kinyon et. al. \cite{phi3},
\cite{phi4}, \cite{phi5}, Ramamurthi and Solarin \cite{ram},
Fenyves \cite{fen1}, \cite{fen2}, Kunen \cite{ken2}, Adeniran
\cite{ade}, Solarin and Chiboka \cite{sol}, \cite{chi1} and
Jaiy\'eol\'a and Ad\'en\'iran \cite{top}. The difficulty in
studying them is as a result of the nature of the identities
defining them when compared with other Bol-Moufang identities. It
can be noticed that in the aforementioned LC identity, the
repeated variable $x$ variables is consecutively positioned and
neither the variable $y$ nor the variable $z$ is between them. A
similarly observation is true in the other two identities(i.e the
RC and C identities). But this observation is not true in the
identities defining Bol loops, Moufang loops and extra loops.
Fenyves \cite{fen2} gave three equivalent identities that define
LC-loops, three equivalent identities that define RC-loops and
only one identity that defines C-loops. But recently, Phillips and
Vojt\v echovsk\'y \cite{phi1}, \cite{phi2} gave four equivalent
identities that define LC-loops and four equivalent identities
that define RC-loops. Three of the four identities given by
Phillips and Vojt\v echovsk\'y are the same as the three already
given by Fenyves. The news ones discovered by Phillips and Vojt\v
echovsk\'y  are as follows.
\begin{equation}\label{lc4:def}
(y\cdot xx)z=y(x\cdot xz)~~\textrm{left central identity}
\end{equation}
\begin{equation}\label{rc4:def}
(yx\cdot x)z=y(xx\cdot z)~~\textrm{right central identity}
\end{equation}
LC-loops, RC-loops, C-loops and their isotopes were first studied in
\cite{top} and it was shown that isotopes exist under triples of the
form $(A,B,B)$ and $(A,B,B)$. Our aim in this study is to
investigate if these two types of triples can give rise to a system
of isotopic LC-loops, RC-loops and C-loops.

Here we treat loops in another way by considering their
derivatives. This concept is well described in \cite{pfl1} and it
was used to answer one of the popular questions raised on G-loops.
If $(G,\cdot )$ is a quasigroup, the operation  $(\cdot )$ is
called a function $F$ on $G$ and we write $x\cdot
y=F(x,y)~and~(G,\cdot )=(G,F)$. If we have several quasigroups on
the same carrier set  $G$, then we write
$(G,F_1)~,~(G,F_2)~,~\ldots $ simply as $F_1,F_2,\ldots $.

\begin{mydef}(Page~79, \cite{pfl1})

Let $G$ be the carrier set of quasigroups $F_1,F_2,F_3,\ldots $.
Then
\begin{enumerate}
\item if $\alpha ~:~(G,F_1)\longrightarrow (G,F_2)$ is an isomorphism, we write $F_2=F_1\alpha $.
\item if $(\alpha ,\beta ,\gamma )~ :~F_1\longrightarrow F_2$ is an isotopism, we write $F_2=F_1(\alpha ,\beta ,\gamma )$.
\item if $(U,V,W)$ is an autotopism of $(G,F)$, then we write $F(U,V,W)=F$.
\end{enumerate}
\end{mydef}

If $(G,F)$ is a non-associative quasigroup, then a fixed element
$a\in G$ determines a new operation $(\circ )$ on $G$ such that
\begin{displaymath}
(a\cdot x)\cdot y=a\cdot (x\circ y).
\end{displaymath}
The operation $(\circ )$ depends entirely on $a\in G$. So we write
$(G,\circ )=F^a$ and call it the left derivative of $F$ w.r.t $a$.
So we have $xL_a\cdot y=(x\circ y)L_a$. So we have isotopism
\begin{displaymath}
(L_a,I,L_a)~:~F^a\longrightarrow F~\textrm{or the
isotopism}~(L_a^{-1},I,L_a^{-1})~:~F\longrightarrow F^a.
\end{displaymath}
So we write $F(L_a^{-1},I,L_a^{-1})=F^a$.

Similarly, we also have what we call right derivative of $F$. Fix
$a\in G$ such that
\begin{displaymath}
x\cdot (y\cdot a)=(x\ast y)\cdot a.
\end{displaymath}
Then we call $(G,\ast )=F_a$ the right derivative of $F$ w.r.t $a$.
So we have
\begin{displaymath}
x\cdot yR_a=(x\ast y)R_a\Rightarrow (I,R_a,R_a)~:~F_a\longrightarrow
F\Rightarrow
\end{displaymath}
\begin{displaymath}
(I,R_a^{-1},R_a^{-1})~:~F\longrightarrow F_a\Rightarrow
F(I,R_a^{-1},R_a^{-1})=F_a.
\end{displaymath}
If $(G,F')$ is a principal isotope of $(G,F)$ of the type
$(R_g,L_f,I)$ , then we call $F'$, $f,g-$isotope of $F$. Thus
$F(R_g,L_f,I)=F'=F_{f,g}$ by notation.\\

In section~3, if a loop $L$ is commutative, it is shown that $L$
been an LC(RC)-loop is necessary and sufficient for it to be an
RC(LC or C)-loop. Consequently, the right(left) derivative,
$a^{-1},e-$ and $e,a^{-1}-$ isotopes of such a C-loop are found to
be C-loops. Furthermore, for such an LC(RC or C)-loop $(L,F)$, it is
shown that
\begin{displaymath}
\{F,F^{a^{-1}},F_{a^{-1},e}\}~\textrm{and}~
\{F,F_{a^{-1}},F_{e,a^{-1}}\}
\end{displaymath}
are systems of isotopic C-loops that obey a form of generalized
distributive law. It is proved that for a loop $(L,\theta )$ to be
an LC(RC or C)-loop, it is necessary and sufficient for the
parastrophe $(L,\theta ^*)$ to be a RC(LC or C)-loop. Hence,
isotopes $(L,\otimes )$ and $(L,\ominus )$ of $(L,\theta )$ and
$(L,\theta ^*)$ respectively are proved to be isotopic if either
$(L,\otimes )$ or $(L,\ominus )$ is commutative. The multiplication
tables of an LC-quasigroup isotope, RC-loop parastrophe of an
LC-loop in \cite{fen2} and an RC-quasigroup isotope of the RC-loop
parastrophe are constructed in this section. And these four were
found to form a system of isotopic C-quasigroup under the above
mentioned condition.

Section~4 generalizes the part of the results in \cite{top} on group
isotopes of C-loops. Under the triple(s) mentioned in \cite{top}, it
is shown that C-loops are isotopic to some finite indecomposable
groups of the classes ${\cal D}_i,i=1,2,3,4,5$ in \cite{orin4}.
Whence, it is proved that the center of such a C-loop has a rank of
1,2 or 3.

For other concepts on Loop Theory mentioned in this paper, readers
should please consult Bruck\cite{b2} and Pflugfelder\cite{pfl1}.

\section{Preliminaries}
The following results which are proved in \cite{top} will be
judiciously used in this study.
\begin{myth}\label{iso:lcrc}
Let $(G, \cdot )$ and $(H, \circ )$ be any two distinct loops. If
the triple $\alpha =(A, B, B)~\Big(\alpha =(A, B, A)\Big)$ is an
isotopism of $(G, \cdot )$ upon $(H, \circ )$, then $(G, \cdot )$ is
an LC(RC)-loop $\Leftrightarrow$ $(H, \circ )$ is a LC(RC)-loop.
\end{myth}

\begin{myth}\label{iso:c}
Let $(G, \cdot )$ and $(H, \circ )$ be two distinct loops. If $G$ is
a central square C-loop, $H$ an alternative central square loop and
the triple $\alpha =(A, B, B)~\Big(\alpha =(A, B, A)\Big)$ is an
isotopism of $G$ upon $H$, then $H$ is a C-loop.
\end{myth}

\begin{myth}\label{iso:cc}
Let $(G, \cdot )$ and $(H, \circ )$ be commutative loops. If $\alpha
=(A, B, B)$ or $\alpha =(A, B, A)$ is an isotopism of $G$ upon $H$,
then $G$ is a C-loop if and only if $H$ is a C-loop.
\end{myth}

\section{Derivatives and Parastrophes}

\begin{mylem}\label{lc:rc}
Let $(L,\cdot )$ be a commutative loop. $L$ is an LC(RC)-loop $\Leftrightarrow $ it is an
RC(LC)-loop.
\end{mylem}
{\bf Proof} \\
Let $L$ be a commutative LC(RC)-loop, then by identity (\ref{lc2:def})((\ref{rc2:def})) :
$(x\cdot xy)z=x(x\cdot yz)\Rightarrow z(x\cdot xy)=(x\cdot yz)x\Rightarrow
z(x\cdot yx)=(x\cdot zy)x\Rightarrow z(yx\cdot x)=(zy\cdot x)x\Leftrightarrow $
$L$ is an RC(LC)-loop by identity (\ref{rc2:def})((\ref{lc2:def})).
Hence the proof is complete.

\begin{mycor}\label{rclc:c}
A commutative loop $L$ is an RC(LC)-loop $\Leftrightarrow $ it is a C-loop.
\end{mycor}
{\bf Proof} \\
Recall that in \cite{fen2}, a loop is a C-loop if and only if it is
both an RC-loop and an LC-loop. Using this fact alone, the
sufficient part is proved. The necessary part follows by this fact
and Lemma~\ref{lc:rc}.

\subsection*{Derivatives of Central Loops}

\begin{myth}\label{lc:der}
Let $(L,F)$ be a commutative LC-loop with identity $e$. If $(L,F_1)$
is a $(I,B,B)$-commutative loop isotope of $(L,F)$ such that $B\in \Pi (F)$, then
\begin{enumerate}
\item $F$ is a commutative C-loop.
\item $F_1=F^{a^{-1}}=F_{a^{-1}}$ is a commutative C-loop for any fixed $a\in L$.
\item $F_1\cong F_{a^{-1},e}\cong F_{e,a^{-1}}$ ; hence $F_{a^{-1},e}$ and $F_{e,a^{-1}}$
are commutative C-loops.
\end{enumerate}
\end{myth}
{\bf Proof} \\
\begin{enumerate}
\item $F$ is a commutative C-loop follows immediately from
Corollary~\ref{rclc:c}.
\item If $(I,B,B)~:~F\to F_1$ is an isotopism,
then $x\circ yB=(x\cdot y)B$. With $B\in \Pi(L)$, there exists $a\in
F$ such that $B=R_a=L_a$, hence
\begin{displaymath}
x\circ yR_a=(x\cdot y)R_a\Rightarrow (I,R_a,R_a)~:~F\longrightarrow
F_1~\textrm{is an isotopism}.
\end{displaymath}
$F$ is an RC-loop $\Rightarrow ~R_a=R_{a^{-1}}^{-1}$. Thus
\begin{displaymath}
(I,R_{a^{-1}}^{-1},R_{a^{-1}}^{-1})~:~F\longrightarrow
F_1\Rightarrow F_1=F(I,R_{a^{-1}}^{-1},R_{a^{-1}}^{-1})=F_{a^{-1}}.
\end{displaymath}
$(B,I,B)~:~F\longrightarrow F_1$ is an isotopism implies
$(L_a,I,L_a)~:~F\longrightarrow F_1$ is an isotopism which implies
\begin{displaymath}
(L_{a^{-1}}^{-1},I,L_{a^{-1}}^{-1})~:~F\longrightarrow
F_1~\textrm{is an isotopism }
\end{displaymath}
\begin{displaymath}
\textrm{which also implies}~
F_1=F(L_{a^{-1}}^{-1},I,L_{a^{-1}}^{-1})\Rightarrow F_1=F^{a^{-1}}.
\end{displaymath}
Whence $F_1=F_{a^{-1}}=F^{a^{-1}}$. By Theorem~\ref{iso:lcrc}, $F_1$
is an RC-loop. By hypothesis, $F_1$ is commutative, hence by
Corollary~\ref{rclc:c}, $F_1$ is a commutative C-loop.
\item Recall that $F_{f,g}=F(R_g,L_f,I)$. Let $f=a^{-1}~,~g=e$, then
\begin{displaymath}
F_{a^{-1}~,e}=F(R_e,L_f,I)=F(I,L_{a^{-1}},I).~\textrm{But}~
F^{a^{-1}}=F(L_{a^{-1}}^{-1},I,L_{a^{-1}}^{-1}),~\textrm{thus}
\end{displaymath}
\begin{displaymath}
F^{a^{-1}}=F(I,L_{a^{-1}},I)L_{a^{-1}}^{-1}=F_{a^{-1}~,e}L_{a^{-1}}^{-1}\Rightarrow
F^{a^{-1}}\cong F_{a^{-1}~,e}.
\end{displaymath}
Similarly, if $f=e~,~g=a^{-1}$, then
\begin{displaymath}
F_{e,a^{-1}}=F(R_{a^{-1}},L_e,I)=F(R_{a^{-1}},I,I).~\textrm{But}~
F_{a^{-1}}=F(I,R_{a^{-1}}^{-1},R_{a^{-1}}^{-1}),~\textrm{thus}
\end{displaymath}
\begin{displaymath}
F_{a^{-1}}=F(R_{a^{-1}},I,I)R_{a^{-1}}^{-1}=F_{e,a^{-1}}R_{a^{-1}}^{-1}\Rightarrow
F_{a^{-1}}\cong F_{e,a^{-1}}.
\end{displaymath}
From (2), $F_1=F^{a^{-1}}=F_{a^{-1}}$, thence $F_1\cong
F_{a^{-1}~e}\cong F_{e,a^{-1}}$. $F_1$ is a commutative C-loop,
therefore $F_{a^{-1}~e}$ and $F_{e,a^{-1}}$ are C-loops.
Commutativity is an isomorphic invariant property, hence we conclude
that $F_{a^{-1}~e}$ and $F_{e,a^{-1}}$ are commutative.
\end{enumerate}

\begin{myrem}This result is also true if $(L,F)$ is an RC-loop.
Initially, when we started our study of isotopic invariance of
central loops in \cite{top}, we were unable to find a convenient
$f,g-$isotope to work with which was not the case for Moufang loops
in \cite{pfl1}. But by Theorem~\ref{lc:der}, it is now clear that
$a^{-1},e-$isotopes or $e,a^{-1}-$isotopes are such $f,g-$isotopes
for some $a\in L$.
\end{myrem}

\begin{myth}\label{c:der}
Let $(L,F)$ be a commutative C-loop with identity $e$. If $(L,F_1)$ is a
$(B,I,B)$-commutative loop isotope of $(L,F)$ such that $B\in \Pi (L)$, then
\begin{enumerate}
\item $F_1=F^{a^{-1}}=F_{a^{-1}}$ is a commutative C-loop for any fixed $a\in L$.
\item $F_1\cong F_{a^{-1},e}\cong F_{e,a^{-1}}$ ; hence $F_{a^{-1},e}$ and $F_{e,a^{-1}}$
are commutative C-loops.
\end{enumerate}
\end{myth}
{\bf Proof} \\
By \cite{fen2}, a C-loop is an LC-loop.
Hence the result follows from Theorem~\ref{lc:der}.

\subsection*{System of Isotopic Central Loops}

\begin{myth}\label{lcrcc:sys}
Let $(L,\cdot )=(L,F)$ be a commutative loop with identity $e$. Let
$(L,\circ )=(L,F_1)$ be a $(I,B,B)$-commutative loop isotope of $(L,F)$ such that $B\in \Pi (L)$.
If $(L,\cdot )$ is an LC-loop, RC-loop or C-loop, then
\begin{displaymath}
\{F~,~F^{a^{-1}}~,~F_{a^{-1},e}\}~\textrm{and}~
\{F~,~F_{a^{-1}}~,~F_{e,a^{-1}}\}
\end{displaymath}
are systems of isotopic commutative C-loops that obey a form of
generalized distributive law.
\end{myth}
{\bf Proof} \\
From the results in Theorem~\ref{lc:der} and Theorem~\ref{c:der},
$F_1=F^{a^{-1}}=F_{a^{-1}}$ is a commutative C-loop for any fixed $a\in L$.
\begin{displaymath}
F_1\cong F_{a^{-1},e}\cong F_{e,a^{-1}}~:~\textrm{hence}~
F_{a^{-1},e}~\textrm{and}~F_{e,a^{-1}}
\end{displaymath}
are commutative C-loops.
\begin{displaymath}
F^{a^{-1}}=F(L_{a^{-1}}^{-1},I,L_{a^{-1}}^{-1}),~\textrm{thus}
\end{displaymath}
\begin{displaymath}
F^{a^{-1}}=F(I,L_{a^{-1}},I)L_{a^{-1}}^{-1}=F_{a^{-1}~,e}L_{a^{-1}}^{-1}\Rightarrow
\end{displaymath}
\begin{displaymath}
F^{a^{-1}}=F_{a^{-1}~,e}L_{{(a^{-1})}^{-1}}=F_{a^{-1}~,e}L_a=F_{a^{-1}~,e}B\Rightarrow
B~:~F_{a^{-1}~,e}\cong F^{a^{-1}}.
\end{displaymath}
Thence; if
\begin{displaymath}
F_{a^{-1}~,e}=(L,\ast ),~(x\ast y)B=xB\circ yB.
\end{displaymath}
But $B=L_a=R_a$, thus
\begin{displaymath}
(x\ast y)L_a=xL_a\circ yL_a\Rightarrow a\cdot (x\ast y)=(a\cdot
x)\circ (y\cdot a)\Rightarrow
\end{displaymath}
the system
\begin{displaymath}
\{F~,~F^{a^{-1}}~,~F_{a^{-1},e}\}
\end{displaymath}
obeys a form of the generalized left distributive law for all fixed
$a\in L$. Also
\begin{displaymath}
(x\ast y)R_a=xR_a\circ yR_a\Rightarrow (x\ast y)\cdot a=(x\cdot
a)\circ (y\cdot a)\Rightarrow \{F~,~F^{a^{-1}}~,~F_{a^{-1},e}\}
\end{displaymath}
obeys a form of the generalized right distributive law for all fixed
$a\in L$. Therefore the system
\begin{displaymath}
\{F~,~F^{a^{-1}}~,~F_{a^{-1},e}\}
\end{displaymath}
obeys a form of the generalized distributive law for all fixed $a\in
L$.

On the other hand,
\begin{displaymath}
F_{a^{-1}}=F(I,R_{a^{-1}}^{-1},R_{a^{-1}}^{-1}),~\textrm{thus}
\end{displaymath}
\begin{displaymath}
F_{a^{-1}}=F(R_{a^{-1}},I,I)R_{a^{-1}}^{-1}=F_{e,a^{-1}}R_{a^{-1}}^{-1}\Rightarrow
F_{a^{-1}}=F_{e,a^{-1}}R_{({a^{-1}})^{-1}}\Rightarrow
\end{displaymath}
\begin{displaymath}
F_{a^{-1}}=F_{e,a^{-1}}R_a\Rightarrow
F_{a^{-1}}=F_{e,a^{-1}}B\Rightarrow B~:~F_{e,a^{-1}}\cong
F_{a^{-1}}.
\end{displaymath}
Thence ; if
\begin{displaymath}
F_{e,a^{-1}}=(L,\star ),~(x\star y)B=xB\circ yB.
\end{displaymath}
But $B=L_a=R_a$, thus
\begin{displaymath}
(x\star y)L_a=xL_a\circ yL_a\Rightarrow a\cdot (x\star y)=(a\cdot
x)\circ (y\cdot a)\Rightarrow
\end{displaymath}
the system
\begin{displaymath}
\{F~,~F_{a^{-1}}~,~F_{e,a^{-1}}\}
\end{displaymath}
obeys a form of the generalized left distributive law for all fixed
$a\in L$. Also
\begin{displaymath}
(x\star y)R_a=xR_a\circ yR_a\Rightarrow (x\star y)\cdot a=(x\cdot
a)\circ (y\cdot a)\Rightarrow \{F~,~F_{a^{-1}}~,~F_{e,a^{-1}}\}
\end{displaymath}
obeys a form of the generalized right distributive law for all fixed
$a\in L$. Therefore the system
\begin{displaymath}
\{F~,~F_{a^{-1}}~,~F_{e,a^{-1}}\}
\end{displaymath}
obeys a form of the generalized distributive law for all fixed $a\in
L$. The proof is complete.

\begin{myth}\label{lc:sys}
Let $(L,F)$ be a LC-loop with identity $e$. If $(L,F_1)$ is a
$(A,I,A)$-loop isotope of $(L,F)$ such that $A\in \Pi_\lambda (L)$, then
\begin{enumerate}
\item $F_1=F^{a^{-1}}$ for a fixed $a\in L$.
\item $F_1\cong F_{a^{-1},e}$.
\item the system $\{F~,~F^{a^{-1}}~,~F_{a^{-1},e}\}$
of isotopic loops obey a form of generalized left distributive law
for all fixed $a\in L$.
\end{enumerate}
\end{myth}
{\bf Proof} \\
\begin{enumerate}
\item If $(A,I,A)~:~F\longrightarrow F_1$ is an isotopism, then $xA\circ y=(x\cdot
y)A$. $A\in \Pi_\lambda (F)\Rightarrow ~\exists~a\in L~\ni~A=L_a$,
thus
\begin{displaymath}
xL_a\circ y=(x\cdot y)L_a\Rightarrow
(L_a,I,L_a)=(L_{a^{-1}}^{-1},I,L_{a^{-1}}^{-1})~:~F\to
F_1\Rightarrow F_1=F(L_{a^{-1}}^{-1},I,L_{a^{-1}}^{-1})=F^{a^{-1}}.
\end{displaymath}
\item Recall that $F_{f,g}=F(R_g,L_f,I)$. Let $g=e~,~f=a^{-1}$~, then
\begin{displaymath}
F_{a^{-1},e}=F(R_e,L_{a^{-1}},I)=F(I,L_{a^{-1}},I).
\end{displaymath}
\begin{displaymath}
F^{a^{-1}}=F(L_{a^{-1}}^{-1},I,L_{a^{-1}}^{-1})=F(I,L_{a^{-1}},I)L_{a^{-1}}^{-1}
=F(I,L_{a^{-1}},I)A=F_{a^{-1},e}A.
\end{displaymath}
Thus, $A~:~F_{a^{-1},e}\cong F^{a^{-1}}$.
\item From (2), if $F_{a^{-1},e}=(L,\ast )$, then
\begin{displaymath}
(x\ast y)A=xA\circ yA\Rightarrow (x\ast y)L_a=xL_a\circ
yL_a\Rightarrow a\cdot (x\ast y)=(a\cdot x)\circ (a\cdot y).
\end{displaymath}
Hence the system
\begin{displaymath}
\{F~,~F^{a^{-1}}~,~F_{a^{-1},e}\}
\end{displaymath}
of isotopic loops obeys a form of generalized left distributive law
for all fixed $a\in L$.
\end{enumerate}

\begin{myth}\label{rc:sys}
Let $(L,F)$ be an RC-loop with identity $e$. If $(L,F_1)$ is an $(I,B,B)$-loop
isotope of $(L,F)$ such that $B\in \Pi_\rho (F)$, then :
\begin{enumerate}
\item $F_1=F_{a^{-1}}$ for a fixed $a\in L$.
\item $F_1\cong F_{e,a^{-1}}$.
\item the system $\{F~,~F_{a^{-1}}~,~F_{e,a^{-1}}\}$
of isotopic loops obey a form of generalized right distributive law
for all fixed $a\in L$.
\end{enumerate}
\end{myth}
{\bf Proof} \\
\begin{enumerate}
\item If $(I,B,B)~:~F\to F_1$ is an isotopism, then $x\circ yB=(x\cdot
y)B$.
\begin{displaymath}
B\in \Pi_\rho (F)\Rightarrow ~\exists~a\in
L~\ni~B=R_a,~\textrm{thus}~x\circ yR_a=(x\cdot y)R_a\Rightarrow
\end{displaymath}
\begin{displaymath}
(I,R_a,R_a)=(I,R_{a^{-1}}^{-1},R_{a^{-1}}^{-1})~:~F\to
F_1\Rightarrow F_1=F(I,R_{a^{-1}}^{-1},R_{a^{-1}}^{-1})=F_{a^{-1}}.
\end{displaymath}
\item Recall that $F_{f,g}=F(R_g,L_f,I)$. Let $g=a^{-1}~,~f=e$~, then
\begin{displaymath}
F_{e,a^{-1}}=F(R_{a^{-1}},L_e,I)=F(R_{a^{-1}},I,I).
\end{displaymath}
\begin{displaymath}
F_{a^{-1}}=F(I,R_{a^{-1}}^{-1},R_{a^{-1}}^{-1})=F(R_{a^{-1}},I,I)R_{a^{-1}}^{-1}
=F(R_{a^{-1}},I,I)B=F_{e,a^{-1}}B.
\end{displaymath}
Thus, $B~:~F_{e,a^{-1}}\cong F_{a^{-1}}$.
\item From (2), if $F_{e,a^{-1}}=(L,\star )$ then
\begin{displaymath}
(x\star y)B=xB\circ yB\Rightarrow (x\star y)R_a=xR_a\circ
yR_a\Rightarrow  (x\star y)\cdot a=(x\cdot a)\circ (y\cdot a).
\end{displaymath}
Hence the system
\begin{displaymath}
\{F~,~F_{a^{-1}}~,~F_{e,a^{-1}}\}
\end{displaymath}
of isotopic loops obeys a form of generalized right distributive law
for all fixed $a\in L$.
\end{enumerate}

\begin{mycor}\label{gen:sys}
Let $(L,F)$ be a central square C-loop with identity $e$. If :
\begin{enumerate}
\item $(L,F_1)$ is a $(A,I,A)$-loop isotope of $(L,F)$ such that $A\in \Pi_\lambda (L)$, then ;
\begin{description}
\item[(i)] $F_1=F^{a^{-1}}$ for a fixed $a\in L$.
\item[(ii)] $F_1\cong F_{a^{-1},e}$.
\item[(iii)] the system $\{F~,~F^{a^{-1}}~,~F_{a^{-1},e}\}$
of isotopic loops obey a form of generalized left distributive law
for all fixed $a\in L$.
\end{description}
\item $(L,F_1)$ is an $(I,B,B)$-loop isotope of $(L,F)$ such that $B\in \Pi_\rho (F)$, then ;
\begin{description}
\item[(i)] $F_1=F_{a^{-1}}$ for a fixed $a\in L$.
\item[(ii)] $F_1\cong F_{e,a^{-1}}$.
\item[(iii)] the system $\{F~,~F_{a^{-1}}~,~F_{e,a^{-1}}\}$
of isotopic loops obey a form of generalized right distributive law
for all fixed $a\in L$.
\end{description}
\item $(L,F_1)$ is a $(A,I,A)$-alternative central square loop isotope of $(L,F)$ such that $A\in \Pi_\lambda (L)$, then ;
\begin{description}
\item[(i)] $F_1$ is a C-loop.
\item[(ii)] $F_1\cong F_{a^{-1},e}$.
\item[(iii)] the system $\{F~,~F^{a^{-1}}~,~F_{a^{-1},e}\}$
of isotopic central square C-loops that obey a form of generalized
left distributive law for all fixed $a\in L$.
\end{description}
\end{enumerate}
\end{mycor}
{\bf Proof} \\
\begin{enumerate}
\item A loop is a C-loop if and only if it is both an LC-loop and an
RC-loop. Thus using Theorem~\ref{lc:sys} and Theorem~\ref{rc:sys},
the proof of the whole of (1) follows immediately.
\item A loop is a C-loop if and only if it is both an LC-loop and an
RC-loop. Thus using Theorem~\ref{lc:sys} and Theorem~\ref{rc:sys},
the proof of the whole of (2) follows immediately.
\item
\begin{description}
\item[(i)] Since $F$ is a central square C-loop, then by Theorem~\ref{iso:c},
$F_1$ is a C-loop.
\item[(ii)] Following Theorem~\ref{lcrcc:sys},
$F_1\cong F_{a^{-1},e}$ for all fixed $a\in L$ implies
$F_{a^{-1},e}$ is a central square C-loop since by (i), $F_1$ is a
central square C-loop.
\item[(iii)] According to Theorem~\ref{lcrcc:sys}, the system of
isotopic loops ; $\{F~,~F^{a^{-1}}~,~F_{a^{-1},e}\}$ obey a form of
generalized left distributive law for all fixed $a\in L$. In fact by
(i) and (ii), it is a system of central square C-loops.
\end{description}
\end{enumerate}

\begin{myrem}
The set $\Pi_\rho$, $\Pi_\lambda $ or $\Pi$ as used here for a loop is not a groupoid.
But in \cite{rit1}, the author gives examples of regular permutation sets that form loops.
\end{myrem}

\subsection*{A Certain Parastrophe of Central Loop}

The definitions of the five parastrophes or conjugates of a loop are
given in \cite{den}. But here, we are concerned with just one of
them and we give the definition below.

\begin{mydef}\label{par:loop}
Let $(L,\theta )$ be a loop. One of the 5 parastrophes of $(L, \theta )$ with
the binary operation $\theta ^*$ defined on the set $L$ is a loop denoted by
\begin{displaymath}
(L,\theta^*)~: ~y\theta ^*x=z\Leftrightarrow x\theta
y=z~\forall~x,y,z\in L.
\end{displaymath}
\end{mydef}

\begin{mylem}\label{lcrc:par}
Let $(L,\theta )$ be a loop. $(L,\theta )$ is an LC(RC)-loop $\Leftrightarrow $ the
parastrophe $(L,\theta ^*)$ is an RC(LC)-loop.
\end{mylem}
{\bf Proof} \\
Putting in mind, (a) in Definition~\ref{par:loop}, $(L,\theta )$ is an LC(RC)-loop
by identity (\ref{lc1:def})((\ref{rc1:def}))
\begin{displaymath}
\Leftrightarrow ~\forall~x,y,z\in L,~(x\theta (x\theta y))\theta
z=(x\theta x)\theta (y\theta z)\Leftrightarrow
\end{displaymath}
\begin{displaymath}
z\theta ^*((y\theta ^*x)\theta ^*x)=(z\theta ^*y)\theta ^*(x\theta
^*x)~\forall~x,y,z\in L\Leftrightarrow (L,\theta ^*)
\end{displaymath}
is an RC(LC)-loop by identity (\ref{rc1:def})((\ref{lc1:def})). For
the other part $(L,\theta ^*)$ is an RC-loop by
\begin{displaymath}
\Leftrightarrow \forall~x,y,z\in L,~z\theta ^*((y\theta ^*x)\theta
^*x)=(z\theta ^*y)\theta ^*(x\theta ^*x)\Leftrightarrow
\end{displaymath}
\begin{displaymath}
(x\theta (x\theta y))\theta z=(x\theta x)\theta (y\theta z)~
\forall~x,y,z\in L \Leftrightarrow (L,\theta )
\end{displaymath}
is an LC-loop. End of the proof.

\begin{mylem}\label{c:par}
Let $(L,\theta )$ be a loop. $(L,\theta )$ is a C-loop $\Leftrightarrow $ the
parastrophe $(L,\theta ^*)$ is a C-loop.
\end{mylem}
{\bf Proof} \\
Putting in mind, (a) in Definition~\ref{par:loop}, $(L,\theta )$ is a C-loop
by identity (\ref{c:def})
\begin{displaymath}
\Leftrightarrow ~\forall~x,y,z\in L,~((y\theta x)\theta x)\theta
z=y\theta (x\theta (x\theta z))\Leftrightarrow
\end{displaymath}
\begin{displaymath}
z\theta ^*(x\theta ^*(x\theta ^*y))=((z\theta ^*x)\theta ^*x)\theta
^*y~\forall~x,y,z\in L \Leftrightarrow (L,\theta ^*)
\end{displaymath}
is an C-loop by identity (\ref{c:def}).

\subsection*{Construction of Two Isotopic Loops}

\begin{myth}\label{iso:par}
Let $(L,\otimes )$ and $(L,\oplus )$ be two distinct quasigroups isotopic to the loops
$(L,\theta )$ and $(L,\theta ^*)$ respectively. If one of $(L,\otimes )$ and $(L,\oplus )$
is commutative, then they are isotopic.
\end{myth}
{\bf Proof} \\ To avoid confusion, let us denote the left
translations of
\begin{displaymath}
(L,\theta )~,~(L,\theta ^*)~,~ (L,\otimes )~,~(L,\oplus )
\end{displaymath}
\begin{displaymath}
\textrm{by}~L_x^\theta ~,~L_x^*~,~L_x^\otimes
~,~L_x^\oplus~\textrm{respectively and the right translations by}
\end{displaymath}
\begin{displaymath}
R_x^\theta ~,~R_x^*~,~R_x^\otimes ~,~R_x^\oplus~\textrm{respectively
for all $x\in L$.}
\end{displaymath}
By (a) in Definition~\ref{par:loop},
\begin{displaymath}
x\theta y=y\theta ^*x\Leftrightarrow L_x^\theta =R_x^*~\forall~x\in
L.
\end{displaymath}
By the general form of the result in \cite{rit}: $(L,\theta )$ is
isotopic to $(L,\otimes )$ under the triple $\alpha
=(A,B,C)\Leftrightarrow $
\begin{equation}\label{eq:l}
L_{xA}^\otimes =B^{-1}L_x^\theta C~\forall~x\in L
\end{equation}
and $(L,\theta ^*)$ is isotopic to $(L,\oplus )$ under the triple
$\beta=(D,E,F)\Leftrightarrow $
\begin{equation}\label{eq:r}
R_{xE}^\oplus =D^{-1}R_x^*F=D^{-1}L_x^\theta F.
\end{equation}
From (\ref{eq:l}), we have
\begin{equation}\label{eq:ln}
L_x^\theta =BL_{xA}^\otimes C^{-1}.
\end{equation}
Hence, combining (\ref{eq:r}) and (\ref{eq:ln});
\begin{displaymath}
R_{xE}^\oplus =D^{-1}L_x^\theta F=D^{-1}(BL_{xA}^\otimes C^{-1})F
=D^{-1}BL_{xA}^\otimes C^{-1}F\Leftrightarrow
\end{displaymath}
\begin{equation}\label{eq:ch}
R_{xA^{-1}E}^\oplus =D^{-1}BL_x^\otimes C^{-1}F\Leftrightarrow
\end{equation}
\begin{displaymath}
R_{xA^{-1}E}^\oplus =D^{-1}BR_x^\otimes C^{-1}F~\ni~(L,\otimes
)~\textrm{is commutative}
\end{displaymath}
\begin{displaymath}
\Leftrightarrow (L,\oplus )~\textrm{is isotopic to}~(L,\otimes )~
\textrm{under the triple}~\gamma
=\Big(B^{-1}D~,~A^{-1}E~,~C^{-1}F\Big)
\end{displaymath}
by the general form of the result in \cite{rit}.

Alternatively by (\ref{eq:ch}),
\begin{displaymath}
R_{xA^{-1}E}^\oplus =D^{-1}BL_x^\otimes C^{-1}F\Leftrightarrow
L_{xA^{-1}E}^\oplus =D^{-1}BL_x^\otimes C^{-1}F~\ni~(L,\oplus )
\end{displaymath}
is commutative
\begin{displaymath}
\Leftrightarrow (L,\oplus )~\textrm{is isotopic to}~(L,\otimes
)~\textrm{under the isotopism}~\mu
=\Big(A^{-1}E~,~B^{-1}D~,~C^{-1}F\Big)
\end{displaymath}
by the general form of the result in \cite{rit}.

\begin{mycor}\label{iso:par1}
Let $(L,\otimes )$ and $(L,\oplus )$ be two distinct quasigroups isotopic to the loops
$(L,\theta )$ and $(L,\theta ^*)$ respectively. If one of $(L,\otimes )$ and $(L,\oplus )$
is commutative, then either of the following is true.
\begin{enumerate}
\item $(L,\theta ^*)$ and $(L,\otimes )$ are isotopic.
\item $(L,\theta )$ and $(L,\oplus )$ are isotopic.
\end{enumerate}
\end{mycor}
{\bf Proof} \\
The proof of (1) and (2) follows from the proof of Theorem~\ref{iso:par}.

\begin{mycor}\label{iso:parcot}
Let $(L,\otimes )$ and $(L,\oplus )$ be two distinct quasigroups
isotopic to the loops $(L,\theta )$ and $(L,\theta ^*)$
respectively. If one of $(L,\otimes )$ and $(L,\oplus )$ is
commutative, then they can be principal isotopes.
\end{mycor}
{\bf Proof} \\
This follows immediately from Theorem~\ref{iso:par}.

\begin{myrem}
We have the following commutative diagrams.
$$ $$
$$(L,\theta )\longrightarrow ^{^{(A,B,C)}} (L,\otimes )$$
$$\searrow ^{^{(E,D,F)}}$$
$$~~~~~~~~~~~~~~~~~~~~~~~~~~~~~~~~~~\downarrow ^{^{(A^{-1}E,B^{-1}D,C^{-1}F)}}$$
$$(L,\theta ^*)\longrightarrow ^{^{(D,E,F)}} (L,\oplus )$$
$$ $$
$$(L,\theta )\longrightarrow ^{^{(A,B,C)}} (L,\otimes )$$
$$~~~~~~~~~~~~~~~~~~~\swarrow ^{^{(B^{-1},A^{-1},C^{-1})}}$$
$$~~~~~~~~~~~~~~~~~~~~~~~~~~~~~~~~~~\downarrow ^{^{(B^{-1}D,A^{-1}E,C^{-1}F)}}$$
$$(L,\theta ^*)\longrightarrow ^{^{(D,E,F)}} (L,\oplus )$$

\end{myrem}

\subsection*{Construction~1} We shall be considering an LC-loop
$(L,\theta )$, of order 6 taken from \cite{fen2}, which is neither a
C-loop nor a Bol loop. By using the left regular permutations of
$(L,\theta )$, i.e the members of $\Pi_\lambda (L,\theta )$, we were
able to calculate the members of $\Pi_\lambda (L,\otimes )$ with the
aid of the general form of the result in \cite{rit}, where
$(L,\otimes )$ is an isotope of $(L,\theta )$ under the isotopism
$\alpha =(A,B,C)$, where $A=(1~5~2~4~3~6)$ and $B=C=(1~3~4~6~5~2)$.
Hence the multiplication table of this quasigroup isotope is shown
below.

\begin{center}
\begin{tabular}{|c||c|c|c|c|c|c|}
\hline
$\otimes $ & 1 & 2 & 3 & 4 & 5 & 6 \\
\hline \hline
1 & 6 & 4 & 5 & 2 & 3 & 1 \\
\hline
2 & 5 & 3 & 2 & 6 & 1 & 4 \\
\hline
3 & 4 & 5 & 6 & 1 & 2 & 3 \\
\hline
4 & 3 & 6 & 1 & 5 & 4 & 4 \\
\hline
5 & 1 & 2 & 3 & 4 & 5 & 6  \\
\hline
6 & 2 & 1 & 4 & 3 & 6 & 5 \\
\hline
\end{tabular}
\end{center}

By Theorem~\ref{iso:lcrc}, $(L,\otimes )$ is an LC-quasigroup with a
left identity element $e_\lambda '=5$ and only obeys the LC-identity
(\ref{lc2:def}) but does not obey other LC-identities
(\ref{lc1:def}), (\ref{lc3:def}) and (\ref{lc4:def}). This can be
verified using computer programs such as GAP, FINITAS, SEM, MACE,
MAGMA, MAGIC, FINDER, PROVER and OTTER. But the package 'LOOPS' is
not designed for investigating Bol-Moufang identities in
quasigroups.

\begin{myrem}\label{fact}
As shown in \cite{ken2}, a quasigroup which obeys (\ref{lc2:def}) is
not an loop. The construction above is an example that confirms this
fact.
\end{myrem}

\subsection*{Construction~2} It is possible to construct an RC-loop
from an LC-loop and vice versa. This is possible by
Lemma~\ref{lcrc:par}. Thus from the LC-loop $(L,\theta )$, we can
construct an RC-loop $(L,\theta ^*)$ of order 6 by following the
statement in \cite{den} which says that the multiplication table of
$(L,\theta ^*)$ is gotten by bordering the resulting latin square
gotten by transposing the latin square of the multiplication table
of $(L,\theta )$.

Hence, the multiplication table for the RC-loop $(L,\theta ^*)$ is :

\begin{center}
\begin{tabular}{|c||c|c|c|c|c|c|}
\hline
$\theta ^*$ & 1 & 2 & 3 & 4 & 5 & 6 \\
\hline \hline
1 & 1 & 2 & 3 & 4 & 5 & 6 \\
\hline
\label{rc:loop}
2 & 2 & 1 & 5 & 3 & 6 & 4 \\
\hline
3 & 3 & 6 & 1 & 2 & 4 & 5 \\
\hline
4 & 4 & 5 & 6 & 1 & 3 & 2 \\
\hline
5 & 5 & 4 & 2 & 6 & 1 & 3 \\
\hline
6 & 6 & 3 & 4 & 5 & 2 & 1 \\
\hline
\end{tabular}
\end{center}

Whence, we go ahead to construct an RC-quasigroup isotope of $(L,\theta ^*)$.
From the elements of $\Pi_\rho (L,\theta ^*)$, we use the general form of the result in \cite{rit},
to compute the elements $\Pi_\rho (L,\oplus )$ of  a quasigroup isotope $(L,\oplus )$ such that
$\alpha =(A,B,C)$ is the isotopism between them, $A=C=(1~5~2~4~3~6)$ and $B=(1~3~4~6~5~2)$.

Thus the multiplication table for the quasigroup isotope $(L,\oplus )$ of $(L,\theta ^*)$ is ;
\begin{center}
\begin{tabular}{|c||c|c|c|c|c|c|}
\hline
$\oplus $ & 1 & 2 & 3 & 4 & 5 & 6 \\
\hline \hline
1 & 6 & 4 & 1 & 3 & 5 & 2 \\
\hline
2 & 3 & 5 & 2 & 4 & 6 & 1 \\
\hline
3 & 2 & 6 & 3 & 1 & 4 & 5 \\
\hline
4 & 5 & 1 & 4 & 2 & 3 & 6 \\
\hline
5 & 4 & 2 & 5 & 6 & 1 & 3 \\
\hline
6 & 1 & 3 & 6 & 5 & 2 & 4 \\
\hline
\end{tabular}
\end{center}

By Theorem~\ref{iso:lcrc}, $(L,\oplus )$ is an RC-quasigroup with a
right identity element $e_\rho '=3$ and obeys the RC-identity
(\ref{rc2:def}) but does not obey other RC-identities
(\ref{rc1:def}), (\ref{rc3:def}) and (\ref{rc4:def}). This can also
be verified using the computer programs mentioned earlier on except
the package 'LOOPS' which is not designed for investigating
Bol-Moufang identities in quasigroups.

\paragraph{}
What could be the relationship between $(L,\otimes )$ and $(L,\oplus
)$? Particularly, are they isotopic? By Theorem~\ref{iso:par}, if
$(L,\otimes )$ or $(L,\oplus )$ is commutative, then they are
isotopic. But $(L,\otimes )$ and $(L,\oplus )$ are LC-quasigroup and
RC-quasigroup respectively. Consequently, if the formal condition is
assumed for both quasigroups, then by Corollary~\ref{lc:rc} and
Corollary~\ref{rclc:c}, $(L,\otimes )$ and $(L,\oplus )$ are
isotopic C-quasigroups. Thus, if $(L,\theta )$ is a C-loop, then by
Lemma~\ref{c:par}, $(L,\theta ^*)$ is a C-loop, whence by
Corollary~\ref{iso:par1},
\begin{displaymath}
\Big\{(L,\theta )~,~(L,\theta ^*)~,~(L,\otimes )~,~(L,\oplus )\Big\}
\end{displaymath}
is a system of isotopic C-quasigroups.

\begin{myrem}
We must note that the two constructions we have considered so far considers
isotopism between loops and quasigroups, which according to \cite{pfl1}
makes sense.
\end{myrem}

\section{Some Finite Indecomposable Group Isotopes of C-loops}

In \cite{orin4}, five classes of groups ${\cal D}_i,i=1,2,3,4,5$
that depend on some positive integer parameters were described.
These groups are used to classify finite indecomposable $RA$-loops.
It was mentioned in \cite{orin4} that groups of the types ${\cal
D}_i,~i=1,2,3,4,5$ are central square.

\begin{mylem}\label{c:cald}
Under a triple of the form $\alpha =(A,B,A)((A,B,B)) $, if an
alternative central square loop $L$ is isotopic to a group of the
type ${\cal D}_i,~i=1,2,3,4,5$, $L$ is a C-loop.
\end{mylem}
{\bf Proof} \\
By the above that groups of the types ${\cal D}_i,~i=1,2,3,4,5$ are
central square and Theorem~\ref{iso:c}, the claim is true.

\begin{mylem}\label{cald:cald}
Under a triple of the form $(A,B,B)((A,B,B))$, if groups of the form
${\cal D}_i,~i=1,2,3,4,5$ are isotopic, the they are C-loops. Hence,
they belong to the same isomorphic class.
\end{mylem}
{\bf Proof} \\
A group is a C-loop. A group in the class ${\cal D}_i,~i=1,2,3,4,5$
is central square, hence by the result in \cite{top} on isotopic
characterization of C-loops, groups of the type ${\cal
D}_i,~i=1,2,3,4,5$ are C-loops. Groups are $G$-loops(\cite{pfl1}),
hence these types of groups are isomorphic.

\begin{myrem}
As stated in \cite{goo2}, $D_4\in {\cal D}_1$ and $Q_8\in {\cal
D}_2$. Hence, Lemma~\ref{c:cald} is a generalization of the facts in
\cite{top} that $D_4$ and $Q_8$ are isotopes of C-loops.
Furthermore, from \cite{top}, we conclude that $D_4\cong Q_8$(which
is also indirectly stated in \cite{serg}). Thence by the latter part
of Lemma~\ref{cald:cald}, Lemma~\ref{cald:cald} is a generalization
of this conclusion. We can pass a remark on Lemma~\ref{cald:cald}
that the five classes of groups that have been mentioned are not
only parametric but also isomorphic.
\end{myrem}

\begin{mydef}(\cite{goo2})
A loop $L$ is indecomposable if it is not the direct product of two proper subloops.

The rank of a finite abelian group A, sometimes written $\rho (A)$, is the
fewest number of generators of $A$.
\end{mydef}

\begin{mycor}\label{c:ig1}
Under a triple of the type $\alpha =(A,B,A)((A,B,B))$, a central
square loop $L$ which is isotopic to an indecomposable finite group
with LC property, a unique non-identity commutator and centre of
rank 1(2)[3] is a C-loop. Whence $\rho \Big(Z(L)\Big)=1(2)[3]$.
\end{mycor}
{\bf Proof} \\
It was shown in \cite{goo2}, that a finite group $D$ is
indecomposable  with LC, a unique non-identity commutator such that
\begin{displaymath}
\rho \Big(Z(D)\Big)=1(2)[3]\Leftrightarrow D\in {\cal
D}_1~\textrm{or}~D\in {\cal D}_2(D\in {\cal D}_3~or~D\in {\cal
D}_4)[D\in {\cal D}_5].
\end{displaymath}
Thus by Lemma~\ref{c:cald}, $L$ is a C-loop.

By \cite{pfl1}, $L$ is isotopic to $D$ implies that
\begin{displaymath}
Z(D)\cong Z(L)\Rightarrow~\textrm{there exists an
isomorphism}~\theta ~:~Z(D)\to Z(L).
\end{displaymath}
If $<d>=Z(D), d\in D$ then there exists $d'\in L~\ni~<d'>=Z(L)$.
Thus $\rho \Big(Z(L)\Big)=1(2)[3]$.

\begin{myrem}
Corollary~\ref{c:ig1} implies that central square C-loops can be
classified by isotopism which according to \cite{pfl1} is
advantageous over isomorphism in the classification of finite loops.
\end{myrem}


\begin{thebibliography}{99}
\bibitem{ade} J. O. Adeniran (2002), {\it The study of properties of
certain class of loops via their Bryant-Schneider group}, Ph.D
thesis, University of Agriculture, Abeokuta.
\bibitem{b2} R.H. Bruck, {\it A survey of Binary Systems},
Springer- Verlag, Berlin 1966.
\bibitem{rit} R. Capodaglio Di Cocco (1993), {\it On Isotopism and Pseudo-Automorphism of the loops}, Bollettino U. M. I. 7, 199--205.
\bibitem{rit1} R. Capodaglio Di Cocco (2003), {\it Regular Permutation Sets and Loops}, Bollettino U. M. I. 8, 617-628.
\bibitem{orin4} O. Chein and E. G. Goodaire (2002), {\it Three-Generator Indecomposable $RA$ loops}, 30, 3559--3564.
\bibitem{chi1} V. O. Chiboka and A. R. T. Solarin (1993), {\it Autotopism characterisation of G-loops}, Scientific Annals of Al.I.Cuza. Univ. , 19--26.
\bibitem{den} J. Dene and A. D. Keedwell (1974), {\it Latin squares and their applications}, The English University press Lts, 549pp.
\bibitem{fen1} F. Fenyves (1968), {\it Extra Loops I}, Publ. Math. Debrecen, 15, 235--238.
\bibitem{fen2} F. Fenyves (1969), {\it Extra Loops II}, Publ. Math. Debrecen, 16, 187--192.
\bibitem{goo2} E. G. Goodaire, E. Jespers and C. P. Milies (1996), {\it Alternative  Loop Rings}, NHMS(184), Elsevier, 387pp.
\bibitem{top} T. G. Jaiy\'eol\'a  and J. O. Ad\'en\'iran {\it On Isotopic Characterization of central loops}(communicated for publication).
\bibitem{ken2} K. Kunen (1996), {\it Quasigroups, Loops and Associative Laws}, J. Alg.185, 194--204.
\bibitem{serg} S. Lang (1993), {\it Algebra}, Addison Wesley publishing company, Inc, 906pp.
\bibitem{pfl1} H. O. Pflugfelder (1990), {\it Quasigroups and Loops : Introduction}, Sigma series in Pure Math. 7, Heldermann Verlag, Berlin, 147pp.
\bibitem{phi4} M. K. Kinyon, K. Kunen, J. D. Phillips (2002), {\it A generalization of Moufang and Steiner loops}, Alg. Univer. 48,1, 81--101.
\bibitem{phi1} J. D. Phillips and P. Vojt\v echovsk\'y (2005), {\it The varieties of loops of Bol-Moufang type}, Alg. Univer. 3(54), 259--383.
\bibitem{phi5}  M. K. Kinyon, J. D. Phillips and P. Vojt\v echovsk\'y (2005), {\it Loops of Bol-Moufang type with a subgroup of index
two}, Bul. Acad. Stiinte Repub. Mold. Mat. 3(49), 71--87.
\bibitem{phi2} J. D. Phillips and P. Vojt\v echovsk\'y (2005), {\it The varieties of quasigroups of Bol-Moufang type : An equational
approach}, J. Alg. 293, 17--33.
\bibitem{phi} J. D. Phillips and P. Vojt\v echovsk\'y (2006), {\it On C-loops}, Publ. Math.
Debrecen. 68, 1-2, 115--137.
\bibitem{phi3} M. K. Kinyon, J. D. Phillips and P. Vojt\v echovsk\'y, {\it C-loops : Extensions and construction}, J. Alg. and Applica. (to appear).
\bibitem{ram} V. S. Ramamurthi and A. R. T. Solarin (1988), {\it On finite right central loops}, Publ. Math. Debrecen, 35, 260--264.
\bibitem{sol} A. R. T. Solarin and V. O. Chiboka (1995), {\it A Note on G-loops}, Collection of Scientific papers of the Faculty of Science, Kragujevac 17, 17--26.
\end{thebibliography}
\end{document}